\documentclass[reqno,12pt,letterpaper]{amsart}
\usepackage{amsmath,amssymb,amsthm,graphicx,mathrsfs,url}
\usepackage[usenames,dvipsnames]{color}
\usepackage[colorlinks=true,linkcolor=Red,citecolor=Green]{hyperref}
\usepackage{amsxtra}

\def\?[#1]{\textbf{[#1]}\marginpar{\Large{\textbf{??}}}}

\let\epsilon=\varepsilon 

\newcommand{\RR}{{\mathbb R}}

\newcommand{\SP }{{\mathbb S}}
\newcommand{\CC}{{\mathbb C}}

\newcommand{\bfe}{\mathbf{e}}

\newtheorem{thm}{Theorem}

\numberwithin{equation}{section}

\DeclareMathOperator{\Div}{div}
\DeclareMathOperator{\Res}{Res}

\let\Im=\Imag
\DeclareMathOperator{\loc}{loc}

\let\Re=\Real

\DeclareMathOperator{\supp}{supp}
\DeclareMathOperator{\vol}{vol}

\DeclareMathOperator{\tr}{tr}

\newcommand{\pa}{\partial}
\newcommand{\la}{\langle}
\newcommand{\ra}{\rangle}
\newcommand{\cL}{\mathcal{L}}

\title[Wave decay for star-shaped obstacles in $ \RR^3 $]{Wave decay for star-shaped obstacles in $ \RR^3 $: papers 
of Morawetz and Ralston revisited}

\author{Peter Hintz}
\email{phintz@berkeley.edu}
\author{Maciej Zworski}
\email{zworski@math.berkeley.edu}
\address{Department of Mathematics, University of California, Berkeley, CA 94720, USA}

\begin{document}

\maketitle

\section{Introduction}

The purpose of this expository note is to revisit Morawetz's method \cite{Mo2} for obtaining a lower bound on the rate of exponential decay of waves for the Dirichlet problem outside star-shaped obstacles, and to discuss the uniqueness of the sphere as the extremizer of Ralston's \cite{Ral} subsequent sharp lower bound.

The bound on the decay rate is essentially the same as lower bound on the distance between {\em scattering resonances}, $ \Res ( \mathcal O ) $, and the real axis (minimal {\em resonance width}) for the Dirichlet Laplacian outside an obstacle $ \mathcal O $. We refer to \cite{dizzy} and \cite{revres} for background, definitions and pointers to the literature. 

Except for \S \ref{s:Hadamard}, our note is  an expanded version of Morawetz's remarkable but not so well known paper \cite{Mo2}. In particular, we want to draw attention to the mysterious inequality \eqref{eq:mor2}. There is a slight change of constants compared to \cite{Mo2}: we were not able to recover the bound \eqref{eq:mor2} with $ 2d $ replaced by $ \frac32 d $ on the right hand side \cite[Theorem 1]{Mo2}.  That results in $ \frac14$ rather than $ \frac13$ in the lower bound on resonance widths \eqref{eq:mor1}.

\begin{thm}
\label{t:1}
Suppose that $ \mathcal O \subset \RR^3 $ is a {\em star-shaped} obstacle
and let $ \Res ( \mathcal O )$ denote the set of scattering poles of 
the Dirichlet realization of $ -\Delta $ on $ \RR^3 \setminus \mathcal O $. 
Then 
\begin{equation}
\label{eq:mor1}
\inf_{ \lambda \in \Res ( \mathcal O ) } | \Im \lambda | > \frac 14 
{\rm{diam}}\, ( 
\mathcal O )^{-1} . 
\end{equation}
\end{thm}

The constant $ \frac14 $ in \eqref{eq:mor1} is {\em far from being optimal}: using the scattering matrix, Ralston \cite{Ral} showed that in {\em any} odd dimension 
\begin{equation}
\label{eq:ral1} 
\inf_{ \lambda \in \Res ( \mathcal O ) } | \Im \lambda | \geq 
2 \, {\rm{diam}}\, ( 
\mathcal O )^{-1} ,
\end{equation}
and this is optimal for the sphere in dimensions three and five -- see below and \S \ref{s:Hadamard}.  For other geometric constants which take energy (that is, $ \Re \lambda $) into account, see Fern\'andez and Lavine \cite{fela}. 

Resonances for the unit sphere in $ \RR^n $  are given by the zeros of 
Hankel functions $ H_{ \ell + \frac n2 - 1 }^{(2)} ( \lambda ) $, 
each with multiplicity given by the dimension of the eigenspace of
$ \ell ( \ell + n-2 ) $ of the spherical Laplacian (thus $ 2 \ell + 1 $ when
$ n =3$). 
When $ n $ is odd, these zeros are given by the zeros of polynomials
$ p_{ \ell + \frac{ n-3}2 }( \lambda )  $ where, 
\[ p_k ( \lambda ) := \sum_{m=0}^{k} \left(\frac{ i} 2 \right)^m 
\frac{ ( m+k)! }{ m! (m-k)! } \lambda^{k-m} ,\]
see \cite[(9.19)]{Tay} and also \cite{Stef}. 

\begin{figure}[htbp]
\begin{center}
\includegraphics[width=0.6\textwidth]{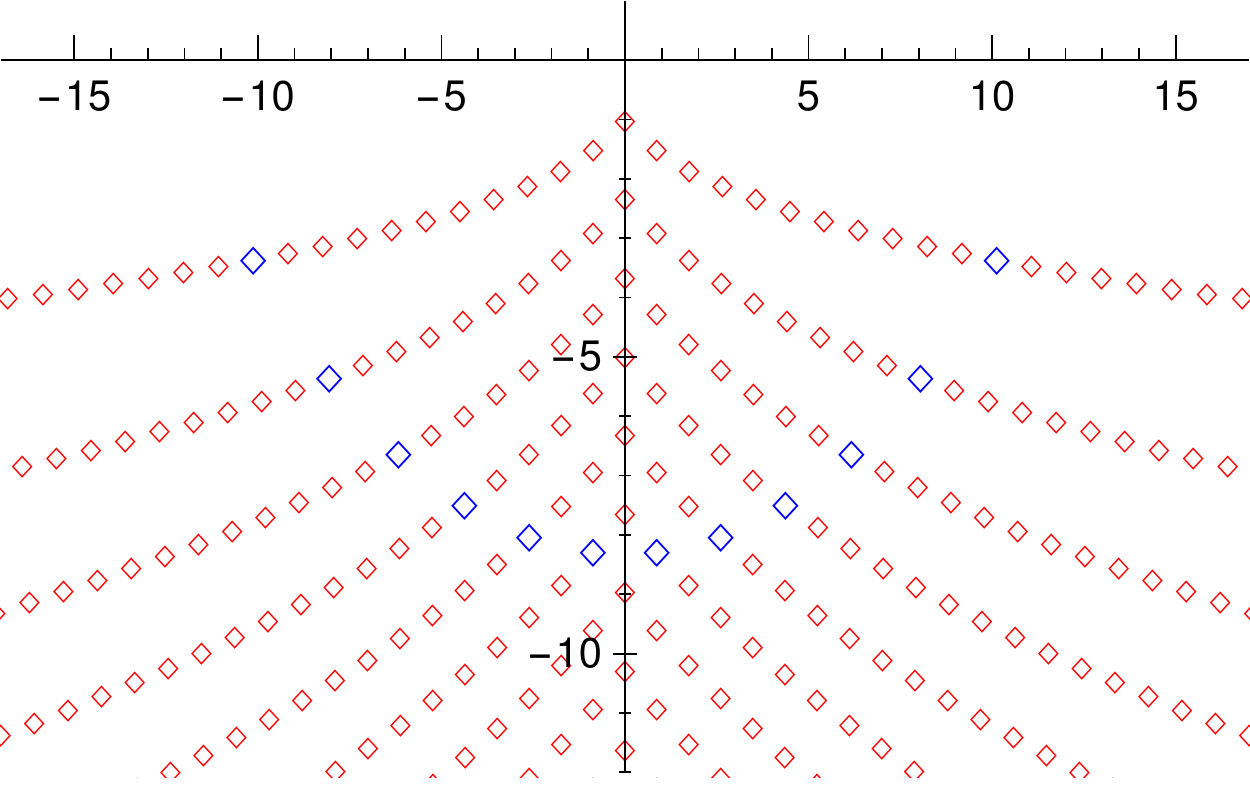}
\caption{\label{f:sphere}
  Resonances for the sphere in three dimensions, see also \cite{Stef}. For each \
spherical
momentum $ \ell $ they are given by solutions of  $ H^{(2)}_{\ell + \frac12} ( \
\lambda ) = 0 $ where $ H^{(2) }_\nu $ is the Hankel function of the second
kind and order $ \nu $. Each zero appears as a resonance of multiplicity
$ 2 \ell +1 $. Highlighted are resonances corresponding to $\ell=20$.}
\end{center}
\end{figure}

One can show (and clearly see from Fig.~\ref{f:sphere}) that for $
n = 3, 5$ the resonance
closest to the real axis comes from solving $ p_1 ( \lambda ) = \lambda + i = 0 $.
That means that 
\begin{equation}
\label{eq:sphere}   \inf_{\lambda \in \Res( B_R ( 0 , 1 ) ) } |\Im \lambda | = 
R^{-1} = 2 \, {\rm{diam}} ( B_R ( 0 , 1 ) )^{-1} , \ \  n = 3,5 ,  
\end{equation}
and Ralston's bound \eqref{eq:ral1} is optimal.

Theorem \ref{t:1} is a consequence of the following theorem, which is valid without
the assumption that $ \mathcal O$ is star-shaped:
\begin{thm}
\label{t:2}
Suppose that $ w $ solves 
\[  ( - \Delta - \lambda^2 ) w = 0 , \ \ x \in \RR^3 \setminus \mathcal O ,\ \ 
w |_{ \partial \mathcal O} = 0,  \]
where $ \mathcal O \subset B ( 0 , d ) $ is an {\em arbitrary} obstacle.

Assume in addition that $ w $ is {\em outgoing} in the sense that 
\begin{equation}
\label{eq:Delw}  w|_{\RR^n \setminus B ( 0 , d )  } = 
( R_0 ( \lambda ) f )|_{\RR^n \setminus B ( 0 , d ) },  \ \ 
\end{equation}
for some $  f \in L^2 ( B ( 0 , d )  ) $, where $R_0(\lambda;x,y)=\frac{e^{i\lambda|x-y|}}{4\pi|x-y|}$ is the integral kernel of the free resolvent. Then 
\[  v ( x ) :=  e^{ - i \lambda |x|} w ( x )  \]
satisfies
\begin{equation}
\label{eq:mor2}
\int_{ \RR^3 \setminus \mathcal O} 
\frac 1 r | \partial_r ( r v ) |^2 dx \leq 2 d
\int_{\RR^3 \setminus \mathcal O} | \partial_x  v |^2 dx .
\end{equation}
\end{thm}
\medskip

\section{Proof of Theorem 1}
We first show how Theorem 2 implies Theorem 1.
For that we first note that 
\[ \begin{split}  e^{ - i \lambda r } \Delta e^{ i \lambda r } & = 
e^{ - i\lambda r } \left( \partial_r^2 + \frac 2 r \partial_r \right) 
e^{i \lambda r }
+ \frac{ \Delta_{ \SP^2 }}{r^2} = 
-\lambda^2 + 2 i \lambda \partial_r + \frac{2 i\lambda}{r} + \Delta  \\
& = - \lambda^2 + \frac{ 2 i \lambda }{r } \partial_r r + \Delta . 
\end{split} \]
Hence, 
if $ ( - \Delta - \lambda^2 ) w = 0 $ in $ \RR^3 \setminus \mathcal O $
and $ w |_{\partial O } = 0 $, then 
\begin{equation}
\label{eq:Delv}   - \Delta v = \frac{2 i \lambda}{r} \partial_r ( r v ) 
\end{equation}
for $ x\in \mathcal E := \RR^3 \setminus \mathcal O $, and $ v |_{\partial \mathcal O } = 0 $.
Multiplying both sides of \eqref{eq:Delv} by $ ( r \bar v )_r $ and taking
real parts we obtain 
\begin{equation}
\label{eq:Delv1}   
\begin{split}
- 2 \Im \lambda \int_{\mathcal E } | ( r v )_r |^2 r^{-1} dx & = - \Re \int_{\mathcal E} \Delta v ( r \bar v )_r d x\\
&  = - \Re \int_{\mathcal E} ( \Delta v \bar v + \Delta v r \partial_r \bar v ) dx \\
& = \int_{ \mathcal E } | \partial_x v |^2 dx +  \int_{\mathcal E} ( - \Re \Delta v r \partial_r \bar v ) dx .
\end{split}
\end{equation}

We put $ F := \partial_x v $ so that the second integrand on the right hand side is
\begin{equation}
\label{eq:divF} \begin{split} - \Re (\partial_x \cdot F ) ( x \cdot \bar F ) & = 
- \Re \partial_x \cdot (F (x \cdot \bar F) ) + \Re F\cdot \partial_x ( x \cdot
\bar F ) \\
& = - \Re \partial_x \cdot (F(x \cdot \bar F) - \tfrac{1}{2}x |F|^2 ) - {\textstyle{\frac12} }|F|^2 .
\end{split} \end{equation}
Here we used the fact that $ F $ is a gradient to obtain the second equality:
\begin{align*}
  \Re \partial_x v \cdot \pa_x( x \cdot \pa_x \bar v ) &= 
\Re \sum_{i,j=1}^3 (\partial_{x_j}v) \pa_{x_j}( x_i \partial_{x_i} \bar v) = 
\sum_{j=1}^3 |\partial_{x_j} v|^2 + {\textstyle{\frac12} } \sum_{ i,j=1}^3
x_i \partial_{x_i} ( | \partial_{x_j} v |^2 ) \\
  &= -\tfrac{1}{2}|\pa_x v|^2 + \tfrac{1}{2}\pa_x\cdot(x|\pa_x v|^2)
\end{align*}
Returning to \eqref{eq:Delv1} and using \eqref{eq:divF} and the divergence theorem, we obtain 
\begin{equation}
\label{eq:Delv2}   
\begin{split}
& - 2 \Im \lambda \int_{\mathcal E } | ( r v )_r |^2 r^{-1} dx \\
& \ \ \ \ =  
{\textstyle \frac12} \int_{\mathcal E} | \partial_x v |^2 dx + \Re \int_{\partial \mathcal E }
(n \cdot \partial_x v  ) ( x \cdot \partial_x v ) d \sigma - 
{\textstyle \frac12} \int_{ \partial \mathcal E } ( x \cdot n ) |\partial_x v|^2 d \sigma ,
\end{split}
\end{equation}
where $ n $ is the {\em outward} (as far as $ \mathcal O $ goes) pointing unit normal vector on $ \partial \mathcal E $ (that is inward pointing for $ \mathcal E $ --- hence
the change of sign). Since $ v |_{\partial \mathcal E } = 0 $, we have $ \partial_x v = n \partial_\nu v $, where the normal derivative is defined by $  \partial_\nu v := n \cdot \partial_x v $; this shows that 
\begin{equation}
\label{eq:Delv3}   
\begin{split}
& - 2 \Im \lambda \int_{\mathcal E } | ( r v )_r |^2 r^{-1} dx =
{\textstyle \frac12} \int_{\mathcal E} | \partial_x v |^2 dx + 
{\textstyle \frac12} \int_{ \partial \mathcal E } ( x \cdot n ) |\partial_\nu v|^2 d \sigma ,
\end{split}
\end{equation}
From Theorem \ref{t:2} we obtain (assuming, as we may, that $ \Im \lambda < 0 $),
\[ 2 | \Im \lambda | \int_{\mathcal E } | ( r v )_r |^2 r^{-1} dx  
\leq 2  | \Im \lambda | {\rm{diam}} ( \mathcal O) \int_{\mathcal E} 
| \partial_x v |^2 dx , \]
which combined with \eqref{eq:Delv3} gives
\[ {\textstyle \frac12} \int_{ \partial \mathcal E } ( x \cdot n ) |\partial_x v|^2 d \sigma \leq {\textstyle{\frac{1}2}} ( 4| \Im \lambda | {\rm{diam}} ( \mathcal O) -1 ) 
\int_{\mathcal E} | \partial_x v |^2 dx  . \]
For a star-shaped obstacle we can choose the origin so that $ x \cdot n \geq 0 $
and hence the left hand side is positive. This gives \eqref{eq:mor1}.

\section{The key estimate}

Suppose that 
\begin{equation}
\label{eq:assu} \Box u ( t, x )  = 0 ,  \  ( t , x ) \in  [ 0 , \infty ) \times \RR^3 , \ \ 
 u (t,  x ) = 0  , \ \ |x | < t - 2 d . 
 \end{equation}
Then
\begin{equation}
\label{eq:mor3}
\begin{split} 
& \Re \int_{ t = d , r \leq d } ( ru_r + u ) \bar u_t dx 
+ {\textstyle{\frac1 { \sqrt 2 } } } \int_{ r = t, t \geq d } 
\left( t | u_t + u_r|^2 + \Re ( u_t + u_r ) \bar u \right) d \sigma 
\\
& \ \ \ \ \leq 
{\textstyle{\frac{1}{2}}} d \int_{t=d, r \leq d } ( |u_x|^2 + |u_t|^2 ) dx 
+ d \int_{ r=t} | \partial_* u |^2 d \sigma - \liminf_{ T \to \infty} 
{ \int_{ r = t = T } | u|^2 dS },
\end{split}
\end{equation}
where $ | \partial_* u |$ denotes the norm of the surface gradient.
This inequality assumes bounds needed to obtain \eqref{eq:limT} below.
These bounds are certainly satisfied in the case of 
$ u ( t, x ) = e^{ i\lambda ( |x| - t ) } v
(x )  $, $ |x| > d $ which will be the case to which \eqref{eq:mor3} is 
applied.

\begin{proof}[Proof of \eqref{eq:mor3}]
We start with the following {\em energy identity} (attributed to Protter in \cite{Mo2}): if
\[  V := x\partial_x + t \partial_t , \]
then
\begin{equation}
\label{eq:prot}
\begin{split}
-\Re \Box \bar u ( Vu  + u ) & = 
\partial_x \cdot \left( - \Re (V u  + u ) \bar u_x + 
\tfrac{1}{2} x ( |u_x|^2 - |u_t|^2 ) \right) \\
& \ \ \ \ + \partial_t \left( \Re ( V u + u ) \bar u_t +
\tfrac{1}{2} t ( |u_x|^2 - |u_t|^2 ) \right),
\end{split}
\end{equation}
where we use the convention $\Box = -\pa_t^2 + \pa_x^2$ -- see \S\ref{SecProtter} for a derivation.

For $ u $ satisfying \eqref{eq:assu} we integrate both sides over 
the region bounded by 
\begin{gather}
\label{eq:defGa}
\begin{gathered}    \Gamma_{d}   \cup \Gamma_{d,T}^+ \cup \Gamma_{d,T}^- , \ \ \ 
\Gamma_d := \{ t = d, \  r \leq d \}, \\ 
\Gamma_{d,T}^+ := \{ r = t, \  d \leq t \leq \tfrac{1}{2}  T \} , \ \ \ 
\Gamma_{d,T}^- := \{ r = T - t , \ \tfrac{1}{2} T \leq t \leq T \},
\end{gathered}
\end{gather}
see Figure~\ref{FigInt}.

\begin{figure}[!ht]
\centering
\includegraphics{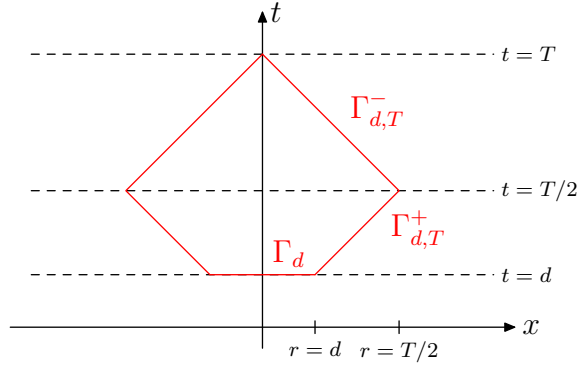}
\caption{Domain of integration.}
\label{FigInt}
\end{figure}

The divergence theorem gives
\begin{equation}
\label{eq:div2}
\begin{split} 
F & = 
\Re \int_{\Gamma_d } ( ru_r + u ) \bar u_t + \tfrac{1}{2} d ( 
| u_x |^2  + |u_t |^2 ) dx \\
& \ \ \ \ + 
{\textstyle{\frac1{\sqrt2}}} \int_{\Gamma_{d,T}^+} \left( 
t | u_t + u_r |^2 + \Re ( u_t + u_r ) \bar u \right) d \sigma , 
\end{split}
\end{equation}
where $ F $ is the contribution from $ \Gamma_{d,T}^- $ (see \eqref{eq:defF}).
The contribution from $ \Gamma_{d,T}^+ $ was calculated as follows: the (Euclidean) outward normal is given by $ ( \mathbf e_r - \mathbf e_t )/\sqrt2 $,
where $ \mathbf e_\bullet $ are the usual unit vectors. Then, since $ r = t $, 
\[ \begin{split} & \mathbf e_r \cdot ( - \Re ( V u + u ) \bar u_x + 
\tfrac{1}{2} x ( |u_x|^2 - |u_t|^2 ) )
- \Re ( V u + u ) \bar u_t -
\tfrac{1}{2} t ( |u_x|^2 - |u_t|^2 ) 
\\ 
& \quad= - \Re ( t u_r + t u_t + u ) \bar u_r  
- \Re ( t u_r + t u_t + u ) \bar u_t 
\\
& \quad= - t | u_t + u_r |^2 - \Re ( u_t + u_r ) \bar u .
\end{split}
\] 
We now calculate the left hand side of \eqref{eq:div2} noting that
the normal vector to $ \Gamma_{d,T}^- $ is $ (\mathbf e_r  +\mathbf e_t )/\sqrt 2 $:
\[ \begin{split} & \mathbf e_r \cdot 
( - \Re ( V u + u ) \bar u_x + 
\tfrac{1}{2} x ( |u_x|^2 - |u_t|^2 ) )
+  \Re ( V u + u ) \bar u_t + 
\tfrac{1}{2} t ( |u_x|^2 - |u_t|^2 ) 
\\ 
& \quad= - \Re ( ( T -  t ) u_r + t u_t + u ) \bar u_r  
+ \Re ( ( T - t )  u_r + t u_t + u ) \bar u_t + \tfrac{1}{2} T ( |u_x|^2 - |u_t|^2 ) 
\\
& \quad= 
t \Re ( |u_r|^2 - 2 u_t \bar u_r +  |u_t|^2 ) 
+ \tfrac{1}{2} T \Re ( - |u_r|^2 + 2 u_r \bar u_t - |u_t|^2 ) 
+ \tfrac{1}{2} T ( |u_x|^2 - |u_r|^2 ) 
\\
&  \quad= 
( t - \tfrac{1}{2} T ) | u_r - u_t|^2 + \Re ( u_t  - u_r ) \bar u 
+ \tfrac{1}{2} T ( |u_x|^2 - |u_r |^2 ) ,
\end{split} \]
so that 
\begin{equation}
\label{eq:defF}
\begin{gathered}
 F = F_1 + F_2 , \ \ F_2 := {\textstyle{\frac1{\sqrt2}}} \int_{\Gamma_{d,T}^-}  
  \Re ( u_t  - u_r ) \bar u  d \sigma , \\
 F_1 := 
{\textstyle{\frac1{\sqrt2}}} \int_{\Gamma_{d,T}^-} 
\left(  ( t - \tfrac{1}{2} T ) | u_r - u_t|^2 
+ \tfrac{1}{2} T ( |u_x|^2 - |u_r |^2 ) \right) d \sigma. 
\end{gathered}
\end{equation}
We start by estimating $ F_2 $: since $ \Re ( u_t - u_r ) \bar u =
\frac12 (\partial_t - \partial_r ) |u|^2 $ 
and $ d \sigma |_{ r = -t + T } = \sqrt 2 ( T - t) ^2 dt d \omega $, $ x = r \omega $, 
we have (recalling that $ u = 0 $ for $ r < T - 2 d $ at $t=T$),
\begin{equation}
\label{eq:F2} \begin{split} F_2 & = {\textstyle{\frac1{2 }}} \int_{t=\frac T 2}^{T} \int_{\SP^2 }
( \partial_t - \partial_r ) |u( t , r\omega ) |^2 |_{ r = T - t } ( T - t)^2 d\omega dt \\
& = - \int_{s=0}^{\frac T 2 } \int_{\SP^2 } \partial_s |u( ( T -s  , s \omega ) |^2 s^2 d\omega ds \\
& =   -\int_{\SP^2 } |u(  {\textstyle{\frac1{2 }}}T, {\textstyle{\frac1{2 }}} T \omega ) |^2 ( {\textstyle{\frac1{2 }}} T )^2 d \omega + 2 \int_{s=0 }^{\frac T 2} \int_{\SP^2 } |u ( T -s  , s \omega )|^2 s d \omega dt \\
& = - \int_{ r = t = \frac12 T } | u|^2 d S +E_T  , 
\end{split}
\end{equation}
where $ dS $ is the surface measure on the sphere defined by 
$ r = t = \frac12 T $ and 
\[ E_T := \sqrt 2  \int_{\Gamma_{d,T}^-} r^{-1} |u|^2 d \sigma .\]

Noting that $ u = 0 $ for $ T - t = r < t - 2 d $ we
see that 
\begin{equation}
\label{eq:supp+}
\Gamma_{d,T}^- \cap \supp u \subset \left\{{\textstyle \frac 12} T  \leq t \leq {\textstyle \frac 12}T + d \right\}
.\end{equation}
Thus, on the support of the integral defining $ E_T $, we have 
$ |r-T/2| \leq d $ and hence
\begin{equation}
\label{eq:ET}
E_T \leq \frac C T \int_{\Gamma_{d,T}^-} |u|^2 d \sigma ,
\end{equation}
and this can be estimated using  \eqref{eq:mor4} and \eqref{eq:Poinc} below.
This shows that 
\begin{equation}
\label{eq:limET}  \lim_{T \to \infty} E_T = 0 .
\end{equation}

We now turn to $F_1 $; using \eqref{eq:supp+} again,
\begin{equation}
\label{eq:F1}
F_1 \leq {\textstyle{\frac1{\sqrt2}}} \int_{\Gamma_{d,T}^-} 
{ {\textstyle{\frac1{2}}} } T ( |  u_x |^2 - | u_r |^2 ) d \sigma
+ {\textstyle{\frac1{\sqrt2}}} \int_{\Gamma_{d,T}^-} 
d | u_r - u_t |^2  d \sigma 
.\end{equation}
Suppose now that $ w $ is another function satisfying \eqref{eq:assu}:
$ \Box w = 0 $ and $ w ( t, x )= 0 $, $ |x| < t - 2d $. We claim that
\begin{equation}
\label{eq:mor4}
 {\textstyle{\frac1{\sqrt2}}} \int_{\Gamma_{d,T}^- } | ( \partial_t - \partial_r ) 
 w |^2 d \sigma \leq 
 \int_{ t = d, r \leq d } (  | w_x|^2 + |w_t|^2 ) dx + 
 {\textstyle{\frac1{\sqrt2}}} \int_{\Gamma_{d,T}^+ } |\partial_* w |^2 d \sigma,
\end{equation}
where $  |\partial_* w | $ is the length of the tangential derivative -- see 
\eqref{eq:deftan}.
For this we use the standard energy identity
\[  - 2 \Re \Box w \bar w_t = \partial_x \cdot ( - 2\Re w_x \bar w_t ) +  
  \partial_t ( |w_x|^2 + |w_t|^2 ) \]
which we integrate over the region bounded by the hypersurfaces in \eqref{eq:defGa}.
That gives (noting that the normals to $ \Gamma_{d,T}^\pm $ are 
$ ( \mathbf e_r \mp \mathbf e_t )/\sqrt 2 $)
\[ \begin{split} 
0 & = -  
\int_{ t = d, r \leq d } (  | w_x|^2 + |w_t|^2 ) dx -
{\textstyle{\frac1{\sqrt2}}}  \int_{\Gamma_{d,T}^+} 
\left( \Re 2 w_r \bar w_t 
 +  ( |w_x|^2 + |w_t|^2 ) \right) 
d \sigma \\
& \ \ \ \ \ + {\textstyle{\frac1{\sqrt2}}} \int_{\Gamma_{d,T}^-} 
\left( - \Re 2 w_r \bar w_t  +  ( |w_x|^2 + |w_t|^2 ) \right)
d \sigma \\
& = -   
\int_{ t = d, r \leq d } (  | w_x|^2 + |w_t|^2 ) dx 
- {\textstyle{\frac1{\sqrt2 }}} \int_{\Gamma_{d , T}^+ } 
\left( | ( \partial_r + \partial_t ) w|^2 + | w_x |^2 - |w_r|^2 \right) d \sigma \\
& \ \ \ \ \ + {\textstyle{\frac1{\sqrt2}}} \int_{\Gamma_{d,T}^-}  
\left(  | ( \partial_r - \partial_t ) w|^2 + | w_x |^2 - |w_r|^2 \right) d \sigma,
\end{split}\]
Since on $ \Gamma_{d,T}^+$, 
\begin{equation}
\label{eq:deftan} |\partial_* w |_{\Gamma_{d,T}^+ }^2 = | ( \partial_r + \partial_t ) w|^2 + | w_x |^2 - |w_r|^2, \end{equation}
we obtain \eqref{eq:mor4}.

We make one more observation: since $ w ( t, ( T - t ) \omega ) $ vanishes for
$  t > \frac12 T + d $ we have
\begin{equation}
\label{eq:Poinc}
\begin{split}
\int_{ \Gamma_{d,T}^- } | w |^2 d \sigma & = 
\sqrt 2 \int_{\frac12 T}^{\frac12 T + d} \int_{\SP^2 } | w ( t , ( T - t ) \omega )|^2 (T-t)^2 d\omega d t \\
& \leq  C_d \int_{\frac12 T}^{\frac12 T + d} \int_{\SP^2 } |   \partial_t w ( t, (T - t)\omega ) |^2 ( T - t)^2 
d \omega dt \\
& = C_d \int_{ \Gamma_{d,T}^- } | ( \partial_t - \partial_r ) w |^2 d \sigma .
\end{split}
\end{equation}
Here we used the following inequality, which holds for $ f $ satisfying $ f ( 0 ) = 0 $
and $ g > 0 $:
\[ \int_0^d |f ( t)|^2 g ( t)  d t = 
\int_0^d \left| \int_0^t f' ( s ) ds \right|^2 g ( t ) dt \leq  
 \frac{\int_0^d g ( t ) t dt }{\min_{t \in [0,d]} g ( t ) } \int_0^d |f'( t ) |^2 g ( t ) dt .\]
(We could compute the $d$-dependent constant but it does not matter as
it disappears in the limit \eqref{eq:limT}.)

We now show that
the first term on the right hand side of \eqref{eq:F1} goes to $ 0 $ as
$ T \to \infty $. To see that we note that on $ \Gamma_{d,T}^+ \cap
\{ 0 \leq t - \frac T2 \leq d \}$,
\[   |u_x |^2 - |u_r|^2 = \frac1 {|x|^2}\sum_{j=1}^3 
| x_j \partial_{x_{j+1} } u  - x_{j+1} \partial_{ x_j} u  |^2 , \  \ x_4 := x_1 ,
\ \ \partial_{x_4 } := \partial_{x_1} .
\] 
Since the vector fields $  x_j \partial_{x_{j+1} }   - x_{j+1} \partial_{ x_j} $
commute with $ \Box $, 
\[ w_j := x_j \partial_{x_{j+1} } u  - x_{j+1} \partial_{ x_j} u, \ \ j = 1, 2, 3,   \]
solve $ \Box w_j = 0 $ and has the same support properties as 
$ u $. Hence to estimate the first term in \eqref{eq:F1}
we can use the estimates \eqref{eq:mor4} and \eqref{eq:Poinc} 
with $ w = w_j$, noting that
on $ \Gamma_{d,T}^- \cap \supp u $, $ |x| \sim T $:
\begin{equation}
\label{eq:limT} \begin{split}  \int_{\Gamma_{d,T}^-} 
  T ( |  u_x |^2 - | u_r |^2 ) d \sigma
& \leq 
\frac{C_d}T \sum_{j=1}^3 \int_{\Gamma_{d,T}^-}  
\sum_{ j=1}^3 | w_j |^2 d \sigma \\
& \leq \frac{C'_d }T \sum_{j=1}^3 \int_{\Gamma_{d,T}^-} 
| ( \partial_t - \partial_r ) w_j |^2 d \sigma \\
& \leq \frac{C'_d }T\sum_{j=1}^3 \left( 
\int_{\Gamma_d}  \sqrt{2}(  | \partial_x w_j|^2 + |\partial_t w_j |^2 ) dx + 
 \int_{\Gamma_{d,T}^+ } |\partial_* w_j |^2 d \sigma
 \right) \\
& \longrightarrow 0 , \ \ T \to \infty.
\end{split}
\end{equation}
Combining this with \eqref{eq:F2}, \eqref{eq:limET}, \eqref{eq:F1} and
using \eqref{eq:mor4} (with $ w = u $) to estimate the {\em second} 
term on the right hand side of \eqref{eq:F1}, we obtain \eqref{eq:mor3}.
\end{proof}

\section{Proof of Theorem \ref{t:2}}

We first show that if 
\begin{equation}
\label{eq:R0f} u_0 := R_0 ( \lambda ) f , \ 
f \in \mathcal D' ( \RR^3 ) , \ \ \supp f \in B ( 0 , d ) , \ \ \lambda \in 
\CC ,
\end{equation}
then the solution of 
\begin{equation}
\label{eq:wavef}
\Box u = 0 , \ \  u |_{t = 0 } = u_0 , \ \ \partial_t u |_{t = 0 } = - i \lambda 
u_0 \end{equation}
satisfies
\begin{equation}
\label{eq:outf}
\supp u \subset \{ ( t, x ) : t < |x| + d  \}.
\end{equation}
This ties the stationary definition of outgoing functions to the dynamical one.

\begin{proof}[Proof of \eqref{eq:outf}]
The argument works of course for any odd $ n \geq 3 $.
We first note that for a fixed $ f $,
 $ \lambda \mapsto u \in C ( \RR_t ; \mathcal D' ( \RR^n ) ) $ is a
holomorphic function. Hence it is enough to prove \eqref{eq:outf} for  
$ \Im \lambda > 0 $ in which case 
$ \hat u_0 ( \xi ) = ( |\xi|^2 - \lambda^2 )^{-1} \hat f ( \xi ) $. Then 
\[ \begin{split} u ( t , x ) & = 
 \left( \cos t \sqrt { - \Delta } - i \lambda \frac {\sin t \sqrt { - \Delta } }
 { \sqrt { - \Delta } } \right) u_0 
\\ &  = \frac{ 1 } { (2 \pi)^n } \int_{ \RR^n } e^{ i \la x , \xi\ra } 
 \left( \cos t |\xi| - i \lambda \frac {\sin t |\xi| }
 { |\xi| } \right) \frac{ \hat f ( \xi ) } { |\xi|^2 - \lambda^2 } d \xi , 
 \end{split} \]
where the Fourier transform is meant in the sense of distributions (the integration 
makes sense for more regular $ f$'s). We can now take the Fourier transform in 
$ t $ which gives, for $ \tau \in \RR $, 
\[  \begin{split} \mathcal F  u ( \tau , x ) & = 
\frac{ 1 } {(2 \pi)^n } \int_{ \RR^n } \int_\RR e^{ i \langle x , \xi \rangle } 
 e^{ - i \tau t } \left( \cos t |\xi| - i \lambda \frac {\sin t |\xi| }
 { |\xi| } \right) \frac{ \hat f ( \xi ) } { |\xi|^2 - \lambda^2 } d t\, d \xi \\
 &  = 
 \frac{ 1 } { 2 (2 \pi)^{n -1} } \sum_{\pm} 
 \int_{ \RR^n } e^{ i \langle \xi, x \rangle } 
 \frac{ \hat f ( \xi ) } { |\xi|^2 - \lambda^2 } 
  \delta ( \tau \mp |\xi| ) ( 1 \mp \lambda / |\xi|) 
  d \xi \\
& = \frac{ 1 } { 2 (2 \pi)^{n -1} } \sum_{\pm} 
 \int_{ \SP^{n-1} }
  e^{ \pm i \tau \langle \omega, x \rangle} 
  \frac{ \hat f ( \pm \tau \omega) }   { \tau^2 - \lambda^2 } 
 ( 1 - \lambda /\tau ) ( \pm \tau)_+^{n-1} 
 d \omega  
\\
& = \frac{ 1 } { 2 (2 \pi)^{n -1} }
\int_{ \SP^{n-1} }
  e^{ i \tau \langle \omega, x \rangle} \frac{ \hat f ( \tau \omega) } { \tau +  \lambda } \tau^{n-2} 
 d \omega  ,
 \end{split} \]
where to get the last equality we crucially used the fact that $ n-1 $ is even.
The expression for $ \mathcal F u ( \tau , x ) $ shows that 
$ \tau \mapsto \mathcal F u ( \tau , x )  $ is holomorphic for $ \Im \tau > - \Im \lambda $
and that, using the Paley--Wiener theorem for $ f $, 
\[ | \mathcal F u ( \tau , x ) | \leq C 
\langle \tau \rangle^M e^{ \Im \tau ( |x| + d ) } . \]
But then \eqref{eq:outf} follows from the Paley--Wiener theorem. 
\end{proof}

Suppose now that $ w $ satisfies the assumptions of Theorem \ref{t:2}, in particular $w=R_0(\lambda)f$ outside of $B(0,d)$, and that 
$ v ( x ) := e^{ - i \lambda |x| } w ( x ) $. Let $ u_0 $ be as in \eqref{eq:R0f}, with the same $f$.
If we solve the free wave equation
\[  \Box U = 0 , \ \ U|_{t=d} = e^{-i\lambda d}u_0 , \ \partial_t U|_{t=d} = - i \lambda e^{-i\lambda d}u_0, \]
then \eqref{eq:outf} shows that $ U $ vanishes for $ |x| < t - 2 d $. Since $ e^{ i \lambda ( |x|
- t ) } v ( x ) $ solves the wave equation in $ \RR \times \{ |x|>d \} $ and
it has the same initial data (at time $ t = d $) as $ U $ in $ |x | > d $ we conclude that
\begin{equation}
\label{eq:u0U}
U ( t , x ) = e^{ i \lambda ( |x| - t) } v ( x) , \ \ |x| \geq t, \ \ t \geq d,
\end{equation}
by the finite speed of propagation property of solutions of the wave equation. Finally we solve the free wave equation $ \Box u = 0  $ with initial conditions
\begin{equation}
\label{eq:solu}
\begin{split}
&   u|_{ t = d } = \begin{cases}  
e^{ i \lambda (|x| -d ) } v ( x ), & |x| > d, \\
v( x ), & x \in \mathcal E \cap \{ |x| \leq d \}, \\
0, & x \in \overline{\mathcal O} ,\end{cases} \\ 
& \partial_t u|_{ t = d } = \begin{cases}  - i \lambda e^{ i \lambda ( |x| -d )} v ( x ),  & |x| 
> d  , \\
0, & |x| \leq d. \end{cases} 
\end{split} \end{equation}
Since $ w|_{\partial \mathcal O} = 0 $, we have $ u |_{ t = d } \in H^1_{\loc} ( \RR^3 ) $, $ \partial_t u  |_{t = d } \in L^2_{\rm{loc}} ( \RR^3 ) $. 

We now apply \eqref{eq:mor3} to $ u ( t,x )$. Since $ u_t |_{t=d} = 0 $ for
$ |x| \leq d $ the first term on the 
left hand side of \eqref{eq:mor3} vanishes.  In the second term $ 
 u ( t, x ) = e^{ i \lambda ( |x| - t )} v ( x ) $ and $ d \sigma = \sqrt 2 dx $.
Hence the left hand side of \eqref{eq:mor3} is given by 
\begin{equation}
\label{eq:lhs1}
\begin{split} 
L & = {\textstyle{\frac1 { \sqrt 2 } } } \int_{ r = t, t \geq d } 
\left( t | u_t + u_r|^2 + \Re ( u_t + u_r ) \bar u \right) d \sigma 
= \int_{ r > d } ( r | v_r |^2 + \Re v_r \bar v ) dx \\
& = \int_{ r > d } \left( r^{-1 } | ( r v)_r |^2 - r^{-1}|v|^2 -
\Re v_r \bar v \right) dx \\
& = \int_{ r > d } \left( r^{-1 } | ( r v)_r |^2 - r^{-1}|v|^2 - 
\tfrac{1}{2}  \partial_r  |v|^2 \right) dx \\
& = 
\int_{ r > d } r^{-1 } | ( r v)_r |^2 dx  -
 \int_{\SP^2}  \int_d^\infty \tfrac{1}{2}  
\partial_r(|v|^2 r^2)  dr d \omega \\
& = \int_{ r > d } r^{-1 } | ( r v)_r |^2 dx + \tfrac{1}{2}   
\int_{ r = d} |v|^2 dS - \lim_{ R \to \infty}  \tfrac{1}{2}   
\int_{ r = R} |v|^2 dS . 
\end{split}
\end{equation}
The right hand side of \eqref{eq:mor3} is
\begin{equation*}
R  =  {\textstyle{\frac{1}{2}}} d \int_{t=d, r \leq d } ( |u_x|^2 + |u_t|^2 ) dx 
+ d \int_{ r=t} | \partial_* u |^2 d \sigma - \liminf_{ T \to \infty} 
{ \int_{ r = t = T } | u|^2 dS } 
\end{equation*}
In view of \eqref{eq:solu} this is equal to 
\begin{equation}
\label{eq:rhs1}
R = {\textstyle{\frac{1}{2}}} d \int_{ \mathcal E \cap  \{ r \leq d \} } | v_x|^2 dx +  d \int_{ \mathcal E \cap \{ r \geq d \} } | v_x|^2 - \lim_{T \to
\infty} \int_{ r = T }  |v|^2 dS .
\end{equation}
Since \eqref{eq:mor3} is $ L \leq  R$ we obtain
\begin{equation}
\label{eq:mor5}
\begin{split} 
& \int_{ r > d } r^{-1 } | ( r v)_r |^2 dx + \tfrac{1}{2}   
\int_{ r = d} |v|^2 dS  \\
& \ \ \ \leq
{\textstyle{\frac{1}{2}}} d \int_{ \mathcal E \cap  \{ r \leq d \} } | v_x|^2 dx + 
 d \int_{ \mathcal E \cap \{ r \geq d \} } | v_x|^2 dx - { {\textstyle{\frac{1}{2}}}  \lim_{R \to
\infty} \int_{ r = R }  |v|^2 dS }\\
& \ \ \ \leq
{\textstyle{\frac{1}{2}}} d \int_{ \mathcal E \cap  \{ r \leq d \} } | v_x|^2 dx + 
 d \int_{ \mathcal E \cap \{ r \geq d \} } | v_x|^2 dx .
 \end{split}
\end{equation}
On the other hand (by integration by parts similar to what we saw before)
\begin{equation}
\begin{split}
\int_{ \mathcal E \cap \{ r < d \} } r^{-1} | ( r v )_r |^2 dx & = 
\int_{ \mathcal E \cap \{ r < d \} } \left( r |v_r|^2 + 2 \Re v_r \bar v + 
r^{-1} |v|^2 \right)dx \\
& = \int_{ \mathcal E \cap \{ r < d \} }  r |v_r|^2 dx + 
\int_{r=d} |v|^2 dS - \int_{ \mathcal E \cap \{ r < d \} } r^{-1} |v|^2 dx \\
& {\leq d \int_{ \mathcal E \cap \{ r < d \} }  |v_x|^2 dx + \int_{r=d} |v|^2 dS}.
\end{split}
\end{equation}
Adding $\frac{1}{2}$ times this inequality to the inequality \eqref{eq:mor5}, we obtain
\[
  \tfrac{1}{2}\int_{\mathcal E} r^{-1}|(r v)_r|^2\,dx \leq d\int_{\mathcal E} |v_x|^2\,dx,
\]
which implies \eqref{eq:mor2}.

\section{Protter's identity from a modern point of view}
\label{SecProtter}

We now explain Protter's identity \eqref{eq:prot} from the point of view presented by Dafermos and Rodnianski \cite[\S4.1.1]{DafermosRodnianskiLectureNotes}, see also \cite{DyatlovQNMExtended}. For that we put
\[
  g := -dt^2 + dx^2.
\]
For $u=u(t,x)$,
\[
  \nabla u = -\pa_t u\,\bfe_t + \nabla_x u,
\]
and for a vector field $V=V_t\,\bfe_t+V_x$ (with $V_x(t,x)$ tangent to $t=t_0$),
\[
  \Div V = \pa_t V_t + \Div_x V_x.
\]
For two vector fields $X$ and $Y$, we introduce
\[
  T_{\nabla u}(X,Y) = \Re\left[g(X,\nabla u)g(Y,\nabla\bar u) - \frac{1}{2}g(X,Y)g(\nabla u,\nabla\bar u)\right].
\]
This defines a new vector field $J_X(u)$ with coefficients quadratic in $\nabla u$ by
\[
  g(J^X(u), Y) = T_{\nabla u}(X,Y).
\]
If $w=w(t,x)$ is a scalar function, one can more generally consider the \emph{modified current}
\[
  J^{X,w}(u) = J^X(u) + \frac{1}{2}\bigl(w\nabla|u^2| - |u|^2\nabla w\bigr),
\]
see for example \cite[\S4.1]{SchlueSchwarzschild}. We then have the following general identity:
\begin{equation}
\label{EqDivId}
\begin{split}
  \Div J^{X,w}(u) &= \Re ((X+w)u)\Box_g\bar u - \frac{1}{2}|u|^2\Box_g w + \Re K^{X,w}(\nabla u,\nabla\bar u), \\
  K^{X,w} &:= \frac{1}{2}\cL_X g - \frac{1}{4}\tr_g(\cL_X g)g + w g.
\end{split}
\end{equation}
If we take $X$ to be the scaling vector field
\[
  X := t\pa_t + x\pa_x,
\]
then
\begin{align*}
  T_{\nabla u}(X,Y)
    &= \Re\bigg\la
         \biggl(
           \begin{bmatrix}
             |\pa_t u|^2           & \pa_t u(\nabla_x\bar u)^T \\
             \pa_t u\nabla_x\bar u & \nabla_x u\otimes\nabla_x\bar u
           \end{bmatrix} \\
   &\qquad\qquad- \frac{1}{2}(-|\pa_t u|^2 + |\nabla_x u|^2)
             \begin{bmatrix}
               -1 & 0 \\
               0  & 1
             \end{bmatrix}
         \biggr)
         \begin{bmatrix}
           X_t \\ X_x
         \end{bmatrix},
       \begin{bmatrix}
         Y_t \\ Y_x
       \end{bmatrix}
     \bigg\ra \\
  &= g(J^X(u), Y) = \bigg\la \begin{bmatrix}-1 & 0 \\ 0 & 1\end{bmatrix}J^X(u), Y \bigg\ra,
\end{align*}
and hence for $X$ as above
\[
  J^X(u)=\Re
  \begin{bmatrix}
    -t |\pa_t u|^2 - \pa_t u x\cdot\nabla_x\bar u + \frac{1}{2}t(|\pa_t u|^2 - |\nabla_x u|^2) \\
    (x\cdot\nabla_x u)\nabla_x\bar u + t\pa_t u\nabla_x\bar u + \frac{1}{2}x(|\pa_t u|^2 - |\nabla_x u|^2)
  \end{bmatrix}.
\]
To compute $K^X=K^{X,0}$, we note that with $\varphi_s(x,t)=(e^s x,e^s t)$,
\[
  \cL_X g = \pa_s\varphi_s^* g|_{s=0} = 2 g
\]
and hence
\[
  K^X = \frac{1}{2}\cL_X g - \frac{1}{4}\tr_g(\cL_X g)g = g - \frac{1}{4}\tr_g(2 g) = -g.
\]
Therefore, if we choose the modifier $w=1$, then $K^{X,w}\equiv 0$, $\Box_g w=0$, and $J^{X,w}(u)=J^X(u) + \Re u\nabla\bar u$, hence the identity~\eqref{EqDivId} becomes
\[
  \Re ((X+1)u)\Box_g\bar u = \Div J^{X,1}(u),
\]
which is exactly Protter's identity~\eqref{eq:prot}.

\section{The variation of the first resonance of the sphere}
\label{s:Hadamard}

We deform $ B ( 0 ,1 ) \subset \RR^3 $ 
without changing the diameter and see the imaginary part of the first
resonance, $ - i $, decreases. In other words, the sphere locally 
maximizes Ralston's bound \eqref{eq:ral1} among obstacles of fixed diameter. 
This result suggests the following

\medskip
\noindent
{\bf Conjecture.} {\em Suppose that $ \mathcal O \subset \RR^3 $ is a non-trapping
obstacle. Then
\[ \inf_{ \lambda \in \Res ( \mathcal O ) } | \Im \lambda | = 1 , \ 
 \mathcal O \subset B ( 0 , 1) 
\ \Longrightarrow \ \mathcal O = B ( 0, 1 ) . \]}
A resolution of this within the class of, say, convex obstacles would already be interesting. At this stage we are not able to gauge the difficulty of this conjecture.

Complex scaling with large angles \cite{SZ1} justifies the following
approach to the variational problem.
We choose a basis of resonant states corresponding to $ - i $ satisfying the following conditions:
\begin{equation}
\label{eq:delij} \int_{\Gamma_\theta}  u_i ( z ) u_j ( z ) dz = \delta_{ ij} , \ \ 
\theta > \pi/2 . \end{equation}
Here the integral is over the radially deformed contour 
(see \cite[(3.16)]{SZ1}) which starts far from the obstacle. Once $ \theta > \pi /2 $ is large enough the integral is independent of $ \theta $
and we drop $ \Gamma_\theta $. We note that $ - \Delta_\theta $ is symmetric with respect to this quadratic form. 

We put $ h ( r ) := r^{-2} e^{r} (r-1) $, 
the radial component of the resonant state corresponding to the resonance at $ - i $.  As spherical harmonics we
choose $ X_j = x_j|_{\SP^2 } $ or explicitly in spherical coordinates
$ ( \theta, \varphi) $, $ 0 \leq \theta \leq 2 \pi $, $ | \varphi | \leq \pi/2 $, $ X_1 = \sin \varphi $, $ X_2 = \cos \varphi \sin \theta $, $ X_3 = 
\sin \varphi \cos \theta $; thus $\int_{\SP^2} X_j^2\,d\!\vol_{\SP^2}=\frac{4\pi}{3}$. With $ A $ to be determined using 
\eqref{eq:delij} we then put
\[ u_j ( r, \theta, \varphi ) = A h ( r ) X_j ( \theta, \varphi ) .\]
We first note that $ \int u_i u_j dz = 0 $ for $ i \neq j $ since we 
complex scale only in the radial variable and the real valued functions $ X_j $
are orthogonal. Now, the integral of $ ( h ( r ) X_j )^2 $ with respect to $ dz$ over $ \Gamma_\theta$, $ \theta > \pi/2$, is
\[ \begin{split}
& \frac {4\pi}{3} \int h ( r ) r^2 dr =  
\frac{4\pi}{3} \int_1^\infty  r^{-2} e^{2r} (r-1)^2 dr
= \frac{4\pi}{3} ( 2r)^{-1} e^{2r} ( r-2)|_{1}^\infty = \frac {2\pi e^2}{3} .
\end{split} \]
Here we can discard the contribution from infinity as we are evaluating the integral over the rescaled contour on which $ e^{2r } $ decays exponentially.
This gives $ A^{-1}  = \sqrt{2\pi e^2/3} $.

We denote by $ z = \lambda^2 $ the ``quantum resonance,'' hence we are deforming $ z = - 1 $ as a Dirichlet eigenvalue of $ - \Delta_\theta $. Since $ - \Delta_\theta $ is symmetric with respect to the quadratic form in \eqref{eq:delij} we can use Hadamard's formula -- see \cite{gri} for a review and references. 
That shows that the first variation comes from 
eigenvalues of the matrix 
\begin{equation}
\label{eq:Cij}  \begin{split} \mathcal C_{ij} & =    \int_{-\pi/2}^{\pi/2}  \int_{0}^{2\pi}  C ( \theta, \varphi ) \partial_r u_j ( 
1 , \theta, \varphi ) \partial_r u_i ( 1, \theta , \varphi ) \cos \varphi d 
\theta d \varphi \\
& = \frac{3}{2\pi} \int_{\SP^2} C \, X_i X_j \, d \! \vol_{\SP^2} , 
\end{split}
\end{equation}
where $ C ( \theta, \varphi ) $ is the normal variation of the obstacle. (The sign difference compared to the standard formula is due to the fact that we 
are applying the formula to the {\em outside} of the obstacle.) Full justification 
comes from a Grushin reduction for the scaled operator and a perturbation formula 
-- see \cite{ela}. 

If a variation does {\em not} increase the diameter of the obstacle
we can assume that the obstacles stay contained in $ B( 0 ,1 ) $. That corresponds to 
\begin{equation}  
\label{eq:Cnp}  C ( \theta, \varphi ) \leq 0. 
\end{equation}

From \eqref{eq:Cij} and \eqref{eq:Cnp}, we see that
\[ \sum_{ i,j} \mathcal C_{ij} \xi_i \xi_j = 
\frac{3}{2\pi} \int_{\SP^2} C \, \langle X, \xi\rangle^2 \, d \! \vol_{\SP^2} , \ \ X := (X_1, X_2, X_3) , \]
and it follows that $ \mathcal C $ is negative semi-definite; if $ C $ is not 
identically zero, $\mathcal C$ is strictly negative.
We conclude that any deformation of the sphere which does \emph{not} increase the diameter moves the first resonance on the imaginary axis deeper into the 
complex half-plane.

To conclude that no other resonance moves closer to the real axis we need to assume a uniform non-trapping condition. Since a smooth deformation has to preserve
convexity for small values of the deformation parameter, 
\cite{hale} and \cite{SZ4} show that resonances lie outside of cubic curves 
determined by the curvature of the obstacle, with the
constants in \cite[(1.3)]{SZ4} depending smoothly on the obstacle. 
Hence continuity of resonances
in compact sets guarantees that all other resonances are at distance more than one from the real axis.

\section{Comparison of the results}

Ralston's proof of \eqref{eq:ral1} uses certain monotonicity properties of the scattering matrix for star-shaped obstacles $\mathcal O\subset\RR^n$.  His argument also allows for suitable perturbations of the Euclidean metric in $B(0,R)$.

\begin{figure}[!ht]
  \centering
  \includegraphics[width=0.6\textwidth]{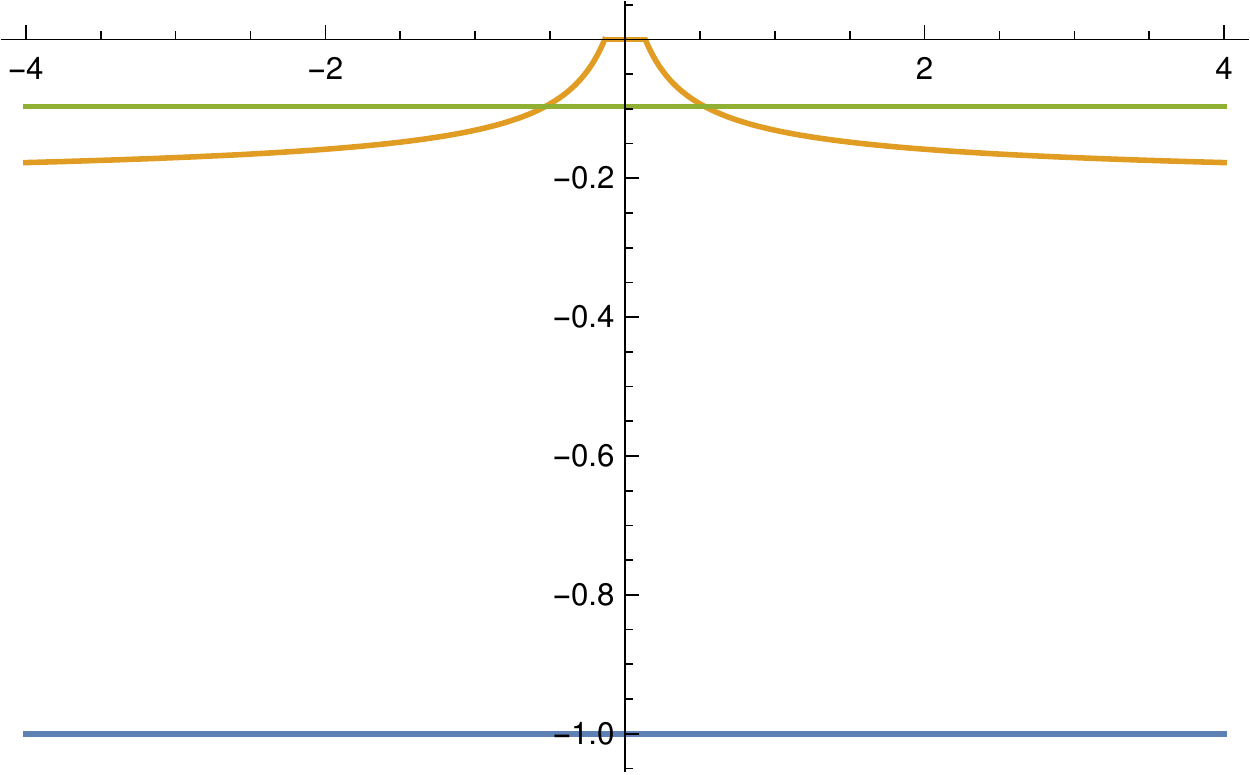}
  \caption{Resonance widths for star-shaped obstacles contained in $B(0,1)\subset\RR^3$ obtained by various authors using various methods. \textit{Blue:} Ralston's unconditional gap. \textit{Yellow:} upper bound for the Fernandez--Lavine gap (setting the $\inf$ in their equation~(5.14) to be equal to $R$). \textit{Green:} the unconditional gap we prove in \cite{HintzZworskiHypobs}.}
\label{FigBounds}
\end{figure}

Fernandez and Lavine \cite{fela} also establish the absence of resonances in certain regions below the real axis, see in particular \cite[Theorem~5.3]{fela} for gaps for obstacle scattering in $\RR^3$ which are however weaker than \eqref{eq:ral1}. Since their methods are different both from those of 
Ralston and Morawetz, we give a brief discussion of their results: due to equation~(5.14) in their paper, their gap becomes worse in particular when the inner radius of the obstacle (the largest ball contained in it) becomes small; the largest possible value of $\alpha$ in (5.14) is thus obtained by replacing the infimum by the constant $R^2$. The bound for $\Im\lambda=:-\eta$ they obtain in their estimate~(5.13) in terms of $\Re\lambda=:\kappa$ is non-trivial unless
\[
  (2\beta\kappa R)^2<3,\quad \beta=1+\frac{e}{2}\Bigl(1+\frac{2}{\kappa R}\Bigr)^{1/2},
\]
which is the case for $\kappa R<0.1353$. As $\Re\lambda\to\infty$, their bound becomes $|\Im\lambda|<\frac{1}{(2+e)R}$, $1/(2+e)\simeq 0.2119$. The different bounds are illustrated in Fig.\ \ref{FigBounds}.

\medskip

\noindent
{\sc Acknowledgements.} MZ would like to thank Cathleen Morawetz for sending him a (rare) reprint of \cite{Mo2} many years ago. PH is grateful to the Miller Institute at the University California, Berkeley for support, and MZ acknowledges partial support under the National Science Foundation grant DMS-1500852. We would also like to thank Jeff Galkowski, Volker Schlue, and Andr\'as Vasy for helpful discussions.  

\def\arXiv#1{\href{http://arxiv.org/abs/#1}{arXiv:#1}}


\begin{thebibliography}{0}

\bibitem[DaRo08]{DafermosRodnianskiLectureNotes}
M.~Dafermos and I.~Rodnianski.
\emph{Lectures on black holes and linear waves,}
in {\em Evolution equations}, Clay Mathematics Proceedings, {\bf 17}(2008),97--205.

\bibitem[Dy11]{DyatlovQNMExtended}
S.~Dyatlov.
\emph{Exponential energy decay for {K}err--de {S}itter black holes beyond
  event horizons,}
\newblock { Mathematical Research Letters}, {\bf 18}(2011), 1023--1035.

\bibitem[DyZw]{dizzy} S.~Dyatlov and M.~Zworski,
		\emph{Mathematical theory of scattering resonances,\/}
                book in preparation; \url{http://math.mit.edu/~dyatlov/res/}

\bibitem[FeLa90]{fela}
C.~Fern\'andez and R.~Lavine,
{\em Lower bounds for resonance widths in potential and obstacle scattering.}
Comm. Math. Phys. {\bf 128}(1990), 263--284.

\bibitem[Gr10]{gri} P.~Grinfeld,
\emph{Hadamard's formula inside and out,}
J. Optim. Theory Appl. {\bf 146}(2010), 654--690.

\bibitem[HaLe94]{hale}
T.~Harg\'e and G.~Lebeau,
{\em Diffraction par un convexe.}
Invent. Math. {\bf 118}(1994),  161--196.

\bibitem[HiZw17]{HintzZworskiHypobs}
P.~Hintz and M.~Zworski,
{\em Resonances for obstacles in hyperbolic space.}
Preprint, 2017.

\bibitem[Mo66a]{Moa} C.~Morawetz, 
\emph{Exponential decay of solutions of the wave equation,} 
Comm.~Pure~Appl.~Math. {\bf 19}(1966), 439--444.


\bibitem[Mo66b]{Moe} C.~Morawetz,
\emph{Energy Identities for the Wave Equation,}
New York Univ., Courant Inst. Math. Sci., Res. Rep. No. IMM346,1966,
\url{https://archive.org/details/energyidentities00mora}


\bibitem[Mo72]{Mo2} C.~Morawetz,
\emph{On the modes of decay for the wave equation in the exterior of 
a reflecting body,} Proc.~Roy.~Irish~Acad.~Sect. A {\bf 72}(1972), 113--120.

\bibitem[Ra78]{Ral} J.~Ralston, 
\emph{Addendum to: "The first variation of the scattering matrix''} (J. Differential Equations {\bf 21}(1976), no.~2, 378--394) by J. W. Helton and Ralston. J. Differential Equations {\bf 28}(1978), no.~1, 155--162. 

\bibitem[Sch13]{SchlueSchwarzschild} V.~Schlue,
\emph{Decay of linear waves on higher-dimensional Schwarzschild black holes,} Anal.~PDE {\bf 6}(2013), no.~3, 515--600.

\bibitem[SjZw91]{SZ1} J. Sj\"ostrand and M. Zworski,
        \emph{Complex scaling and the distribution of scattering poles,\/}
        J. Amer. Math. Soc. \textbf{4}(1991), 729--769.

\bibitem[SjZw95]{SZ4} J. Sj\"ostrand and M. Zworski,
\emph{The complex scaling method for scattering by strictly convex obstacles.}
Ark. Mat. {\bf 33}(1995), 135--172.

\bibitem[SjZw07]{ela} J. Sj\"ostrand and M. Zworski,
        \emph{Elementary linear algebra for advanced spectral problems,\/}
	Ann.~Inst.~Fourier \textbf{57}(2007), 2095--2141.

\bibitem[St06]{Stef} P. Stefanov,
        \emph{Sharp upper bounds on the number of the scattering poles,\/}
        J. Funct. Anal. \textbf{231}(2006), 111--142.

\bibitem[Ta11]{Tay} M.~E.~Taylor, 
{\em Partial Differential Equations II. 
Qualitative Studies of Linear Equations}, Applied Mathematical Sciences
Volume {\bf 116}, Springer, 2011.

\bibitem[Zw17]{revres} M. Zworski,
{\em Mathematical study of scattering resonances,} \arXiv{1609.03550}, to 
appear in Bull.~Math.~Sci.


\end{thebibliography}
\end{document}